\documentclass[11pt]{amsart}
\usepackage{epsfig,amsfonts,amsthm,amssymb,latexsym,amsmath,enumerate,hyperref,color}

\usepackage{epsfig}
\usepackage{epstopdf}

\newtheorem{theorem}{Theorem}

\newtheorem{remark}{Remark}

\title[On an Open Question Concerning Product-Type Difference Equations]
{On an Open Question Concerning Product-Type Difference Equations}

\author[J. F. T. Rabago]{Julius Fergy T. Rabago}
\address{Julius Fergy Tiongson Rabago, Department of Mathematics and Computer Science,
College of Science, University of the Philippines Baguio, Baguio City 2600, Benguet, PHILIPPINES} \email{jfrabago@gmail.com}



\date{\today}
\begin{document}
\maketitle

\begin{abstract}
In [Acta Math. Univ. Comenianae Vol. LXXX, 1 (2011), pp. 63--70], Yang, Chen and Shi examined the system of difference equations
	\[
	x_n=\frac{a}{y_{n-p}},\qquad y_n=\frac{by_{n-p}}{x_{n-q}y_{n-q}},\qquad n=0,1,\ldots,
	 \]    
where $q$ is a positive integer with $p < q$, $p \nmid q$, $p \geq 3$ is an odd number, both $a$ and $b$ are nonzero real constants, and
the initial values $x_{-q+1},x_{-q+2},\ldots,$ $x_0,y_{-q+1},y_{-q+2},\ldots,y_0$ are nonzero real numbers.
At the end of their note, they posted a question regarding the behaviour of solutions of the given system when $p$ is even.
More precisely, they asked what the solutions of the system may look like if $p$ is even.
In this note we answer this question raised by the authors.
Particularly, we show that the system may or may not admit a periodic solution depending on the coprimality of the parameters $p$ and $q$ and on the parity of the integer $p/\gcd(p,q)$. 
\end{abstract}

\section{Introduction}

In a recent paper, Yang, Chen and Shi investigated the system of difference equation
	\begin{equation}
	\label{problem}
	x_n=\frac{a}{y_{n-p}},\qquad y_n=\frac{by_{n-p}}{x_{n-q}y_{n-q}},\qquad n=0,1,\ldots,
	 \end{equation}
where $q$ is a positive integer with $p < q$, $p \nmid q$, $p \geq 3$ is an odd number, both $a$ and $b$ are nonzero real constants, and
the initial values $x_{-q+1},x_{-q+2},\ldots, x_0,$
$y_{-q+1},y_{-q+2},\ldots,y_0$ are nonzero real numbers.
They have showed in \cite{ycs} that all real solutions of the system are eventually periodic with period $2pq$ (resp. $4pq$) when $(a/b)^q = 1$ (resp.
$(a/b)^q = -1$). Further, the authors \cite{ycs} characterized the asymptotic behavior of the solutions of \eqref{problem} when $a \neq b$.
Their work is actually a generalization of a result seen in \cite{oz3}, wherein \"{O}zban investigated the behavior of the positive solutions of the system of rational difference equations
	\begin{equation}
	\label{oz3}
		x_n=\frac{a}{y_{n-3}},\qquad y_n=\frac{by_{n-3}}{x_{n-q}y_{n-q}},\qquad n=0,1,\ldots,
	  \end{equation}
where $q > 3$ is a positive integer with $3\nmid q$, $a$ and $b$ are positive constants.
\"{O}zban particularly showed in \cite{oz3} that for the case $b = a$, $p = 3$, $q > 3$, and $p$ not dividing $q$, 
all positive solutions of the system \eqref{problem} are eventually periodic with period $6q$.
Meanwhile, \"{O}zban's work \cite{oz3} was inspired by a result found by Yang, Liu and Bai in \cite{ylb1}.
In their four page note, Yang et al. \cite{ylb1} investigated the system \eqref{problem} for the case where $p$ and $q$ are positive integers with $p\leq q$, and $a$ and $b$ are positive constants.
They have showed that if $b = a$, and $q$ is divisible by $p$, then every positive solution of \eqref{problem} is periodic with period $2q$.
In \cite{bratislav}, an intriguing result regarding the behavior of positive solutions of the higher-order difference equation 
\begin{equation}
\label{bratislav}
x_n=\frac{cx_{n-p}x_{n-p-q}}{x_{n-q}},\qquad n=0,1,\ldots,
\end{equation} 
where $p,q \in \mathbb{N}$  and $c>0$, was obtained by Iri\v{c}anin and Liu through an elegant and short way.
One particular results they presented in \eqref{bratislav} reads as follows: if $c=1$ in \eqref{bratislav} and $\gcd(p,q)=1$, and $p$ is odd, then all positive solutions of \eqref{bratislav} are eventually periodic with period $2pq$. 
Their result was in fact inspired by earlier results presented in \cite{oz2, oz3} and \cite{ylb1}. 
Similar nonlinear systems of rational difference equations were also investigated, see, e.g., Cinar \cite{c1}, Cinar and Yal\c{c}inkaya \cite{cy1,cy2,cy3}, Cinar et al. \cite{cyr1}, 
and \"{O}zban \cite{oz1}.

Now, our aim in this work is to answer a question raised by Yang et al. at the end of their paper \cite{ycs}.
Particularly, we shall described here the behavior of solutions of \eqref{problem} in the case when $p$ is even.
We mention that Yang et al. \cite{ycs} already made a preliminary observation on the behavior of solutions of \eqref{problem} when $p$ is even.
More specifically, they have observed, after some numerical experimentations, that \eqref{problem} is non-periodic when $p$ is even.
Here we show, through an analytical approach, that this observation is in fact true whenever $\gcd(p,q)=1$, i.e. $q$ is odd
and that, in addition, every positive solution of \eqref{problem} when $b=a$ has an exponentially growing/decaying subsequence irrespective of the choice of positive initial values $x_{-q+1},x_{-q+2},\ldots, x_0,y_{-q+1},y_{-q+2},\ldots,y_0$ on this case.
Every solution, however, will be periodic of period $m$, where $m$ denotes the least common multiple of $p$ and $2q$, if $p/\gcd(p,q)$ is odd.
On the other hand, if $p/\gcd(p,q)$ is even, then every solution of \eqref{problem} behaves in a similar fashion as in the case when $\gcd(p,q)=1$.
That is, there is some subsequence of the solution $\{x_n\}$ (resp. $\{y_n\}$) which tends to infinity (resp. converges to zero) (cf. Theorem \ref{main}).
We emphasize that the problem raised by Yang et al. in \cite{ycs} remains open since Iri\v{c}anin and Liu \cite{bratislav} only deal with the case when $p$ is odd in \eqref{problem}, which is not of interest here.

\begin{remark}
By an eventually periodic solution $\{(x_n,y_n)\}:=\{(x_n,y_n)\}_{n=-(q-1)}^\infty$ of \eqref{problem}, we mean that 
there exist an integer $n_0 \geq -q+1$ and a positive integer $\pi$ such that
\[
(x_{n+n_0+\pi},y_{n+n_0+\pi})=(x_{n+n_0},y_{n+n_0}), \quad n=1,2,\ldots,
\]
and $\pi$ is called a period {\rm (cf. \cite{ladas})}.
If, regardless of the choice of initial values $x_{-q+1}$, $x_{-q+2}, \ldots, x_0, y_{-q+1},y_{-q+2},\ldots,y_0$, no such values of $n_0$ and $\pi$ exist, then every solution of \eqref{problem} is not and can never be periodic. 
\end{remark}
\section{Main Results}

In this section we show that the system \eqref{problem} can never have a periodic or eventually periodic solution when $p$ is even and $q$ is odd.
If, however, $q$ is even, then the existence of periodic solution of \eqref{problem} depends on the parity of $p/\gcd(p,q)$.
Our approach parallels that seen in \cite{bratislav}.

We only consider the case when $b=a$ in \eqref{problem} with all of its initial values taken from the set of positive real numbers. 
The same inductive lines, however, can be followed to show a similar result for the case $b=-a$ and even for the more general case given by system \eqref{problem}.

To show that \eqref{problem} has no periodic solution, it suffices to prove that every solution of \eqref{problem} has an increasing subsequence (or, perhaps, a decreasing subsequence) regardless of the choice of initial values $x_{-q+1},x_{-q+2},\ldots, x_0$, $y_{-q+1},y_{-q+2},\ldots,y_0$. 

With this idea in mind, we now proceed as follows.
First, we transform the first equation in \eqref{problem} to 
\[
x_nx_{n-q}=x_{n-p}x_{n-p-q}, \qquad n=0,1,\ldots. 
\]
Since $x_n>0$ for all $n\geq 0$, then taking the natural logarithm of both sides of the above equation and making the subtitution $a_n:=\ln x_n$, we get
\[
a_n+a_{n-q}-a_{n-p}-a_{n-p-q}=0.
\]
Using the ansatz $a_n =\lambda^n \in \mathbb{R}$, we obtain the polynomial equation
\[
P(\lambda):=\lambda^{p+q}+\lambda^p-\lambda^q-1=(\lambda^p-1)(\lambda^q+1)=0.
\]
From here on, we consider two possibilities: (i) $\gcd(p,q)=1$; and (ii) $\gcd(p,q)>1$.\\

\noindent \underline{CASE 1}: Suppose that $\gcd(p,q)=1$ or equivalently, $q$ is odd. 
Then, it is evident that $\lambda=-1$ is a root of $P(\lambda)=0$ of order two.
Denote this repeated root by $\lambda_1$ and $\lambda_{p+1}$, i.e., let $\lambda_1=\lambda_{p+1}=-1$. 
Then, the explicit formula for the sequence $\{a_n\}$ is of the form
\[
a_n=c_1\lambda_1^n+c_{p+1}n \lambda_{p+1}^n+\sum_{\substack{{i=2}\\{i\neq p+1}}}^{p+q} c_i \lambda_i^n,
\]
for some real numbers $c_1,c_2,\ldots,c_{p+q}$.
Note that $P(\lambda)=0$ has all of its roots on the unit disk $|\zeta|\leq 1$.
Denote these roots by $\{\lambda_i\}_{i=1}^{p+q}$ where $\{\lambda_i\}_{i=1}^p$ are the corresponding roots of $\lambda^p-1=0$  and $\{\lambda_i\}_{i=p+1}^{p+q}$ are the roots of $\lambda^q+1=0$.
Clearly, $\lambda_i^p=1$ for all $1\leq i\leq p$ and $\lambda_i^{2q}=1$ for all $p+1\leq i\leq p+q$.
Since $p \nmid q$, then $\lambda_i^{2pq}=1$ for all $1\leq i\leq p+q$. 
Hence,
	\[
	a_{2pqn}=c_1+2c_{p+1}pqn+\sum_{\substack{{i=2}\\{i\neq p+1}}}^{p+q} c_i.
	\]
Suppose $c_{p+1}>0$, then for sufficiently large $N \in \mathbb{N}$, $a_{2pqn}$ will eventually be increasing for $n\geq N$, in fact we'll have
\[
	a_{2pqn} \longrightarrow \infty \quad \text{as}\quad n \rightarrow \infty.
\]  
Going back to the relation $a_n=\ln x_n$, we see that
\[
x_{2pqn}=\exp\left\{ a_{2pqn}\right\} \longrightarrow \infty \quad \text{as}\quad n \longrightarrow \infty. 
\]
Moreover, we have
 	\[
 	y_{2pqn}\longrightarrow 0\quad \text{as}\quad n \longrightarrow \infty. 
 	\]
Therefore, the subsequence $\{x_{2pqn}\}$ (resp. $\{y_{2pqn}\}$) will tend to infinity (resp. converges to zero) exponentially.
Thus, every solution of \eqref{problem} can never be periodic when $\gcd(p,q)=1$, or equivalently, when $q$ is odd.\\

\noindent\underline{CASE 2}: Now, the case $\gcd(p,q)>1$ needs a little more work. 
Suppose $\gcd(p,q)=2^ru$ for some odd integer $u\geq1$ and integer $r>0$. 
Then, the roots of $P(\lambda)=0$  can be expressed as
\begin{equation}
\label{roots}
\begin{cases}
\ \exp\left\{ \dfrac{(2k+1)\pi i}{q}\right\}, & k=0,1,\ldots, q-1,\\[1em]
\ \exp\left\{ \dfrac{2l \pi i}{p}\right\}, & l=0,1,\ldots, p-1.
\end{cases}
\end{equation}

Let $p=2^rus$ and $q =2^rut$ where $s<t$ and $s\nmid t$.
The roots of $P(\lambda)=0$ are simple if and only if
	\[
		\frac{2k+1}{2^r u t} \neq \frac{2l}{2^r u s}, \qquad \text{for each}\ k,l \in \mathbb{N}_0,
	\]
or equivalently
\begin{equation}
\label{condition}
	(2k+1)s\neq 2l t, \qquad \text{for each}\ k,l \in \mathbb{N}_0.
\end{equation}
We consider two separate subcases, namely: (a) $p/\gcd(p,q)$ is odd; and (b) $p/\gcd(p,q)$ is even.\\
  
\underline{Subcase 2.1}:  Clearly, if $s$ is odd (i.e., $p/\gcd(p,q)$ is odd), then the inequality \eqref{condition} always holds.
Hence, the roots of \eqref{problem} are distinct.
Moreover, since $\lambda^{2q}=1$, then the explicit formula for $a_{mn}$ takes the form     
	\[
	a_{mn+t}=\sum_{i=1}^{p+q} c_i\lambda_i^{mn+t}=\sum_{i=1}^{p+q} c_i\lambda_i^{t}=a_t, \qquad \forall n=0,1,\ldots,
	\]
for each $t=\{0,1,\ldots,m-1\}$, where $m:={\rm lcm}(p,2q)$ denotes the least common multiple (lcm) of $p$ and $2q$. 
Therefore, $a_n$ is eventually periodic with period $m$.
Since $a_n=\ln x_n$ for all $n=0,1,\ldots$, then $x_n$, as well as $y_n$, are also eventually periodic.\\

\underline{Subcase 2.2}: Now, if $p/\gcd(p,q)$ is even (i.e., $s=2^{r_0}s_0<t$ for some integers $r_0,s_0>0$ with $s_0$ being odd), then $t$ must be odd, otherwise $\gcd(p,q) > 2^ru$. 
Evidently, \eqref{condition} does not hold since the equality $(2k+1)2^{r_0-1}s_0= l t$ may hold true by choosing appropriate values for $l,k \in \mathbb{N}_0$.  
For instance, if $r_0=1$, then we can choose $l=s_0$ and $k=(t-1)/2$. 
In general, we can take $l=2^{r_0-1}s_0$ and $k=(t-1)/2$.
Since the inequality was not satisfied, then $P(\lambda)=0$ has at least one repeated root.
Without loss of generality, let $\lambda_j$ be a root of $P(\lambda)=0$ of order two and $m$ be the least common multiple of $p$ and $2q$.
By arguing as in the first case, we see that the subsequence 
\begin{align*}
	a_{mn}&=c_jm n\lambda_j^{mn}+\sum_{i=1}^{p+q} c_i\lambda_i^{mn}=c_j m n+\sum_{i=1}^{p+q} c_i\\
		&\longrightarrow \infty\qquad \text{as}\quad n\longrightarrow \infty.
\end{align*}
Again, going back to the relation $a_n=\ln x_n$, the above result leads us to conclude that when $\gcd(p,q)>1$ and $p/\gcd(p,q)$ is even, 
a subsequence of the solution to \eqref{problem} grows/decays exponentially.
Thus, every solution in this case is non-periodic. \\

We remark that every solution to \eqref{problem} when $b=-a$ has the same behavior as those in the case when $b=a$.
The only difference is that every real solution of \eqref{problem} for $b=-a$ oscillates at $0$.
This can be seen easily from the relation $x_nx_{n-q}=-x_{n-p}x_{n-p-q}$. 
In fact, if $\{(x_n,y_n)\}$ is a solution to \eqref{problem} with positive initial values, then there is some subsequence $\{|x_{mn+t}|\}$ (resp. $\{|y_{mn+t}|\}$)
which tends to infinity (resp. converges to zero) exponentially.\\

We summarize our discussion in the following theorem for the case $b=a$.
A similar conclusion can be established for $b=-a$.

\begin{theorem}
\label{main}
Let $\{(x_n,y_n)\}$ be a solution to \eqref{problem}. 
Then, the the following statements are true:
\begin{enumerate}
\item[(i)] If $\gcd(p,q)=1$, then the solution $\{(x_n,y_n)\}$ of system \eqref{problem} has a subsequence $\{x_{2pqn}\}$ {\rm ({\it resp.} $\{y_{2pqn}\}$)} that tends to infinity {\rm ({\it resp. converges to zero})} exponentially, and vice versa.

\item[(ii)] If $\gcd(p,q)>1$, and $p/\gcd(p,q)$ is odd, then the solution $\{(x_n,y_n)\}$ of system \eqref{problem} is eventually periodic with period $m$, where $m$ denotes the least common multiple of $p$ and $2q$. 

\item[(iii)] If $\gcd(p,q)>1$ and $p/\gcd(p,q)$ is even, then the solution $\{(x_n,y_n)\}$ behaves in a similar fashion as in (i).
\end{enumerate}
\end{theorem}

Some illustrations which shows different behaviors of various solutions of system \eqref{problem} for the case $b=a$ with random positive initial values are illustrated in Figures (1)--(5).
\begin{figure}
    \scalebox{.6}{\includegraphics{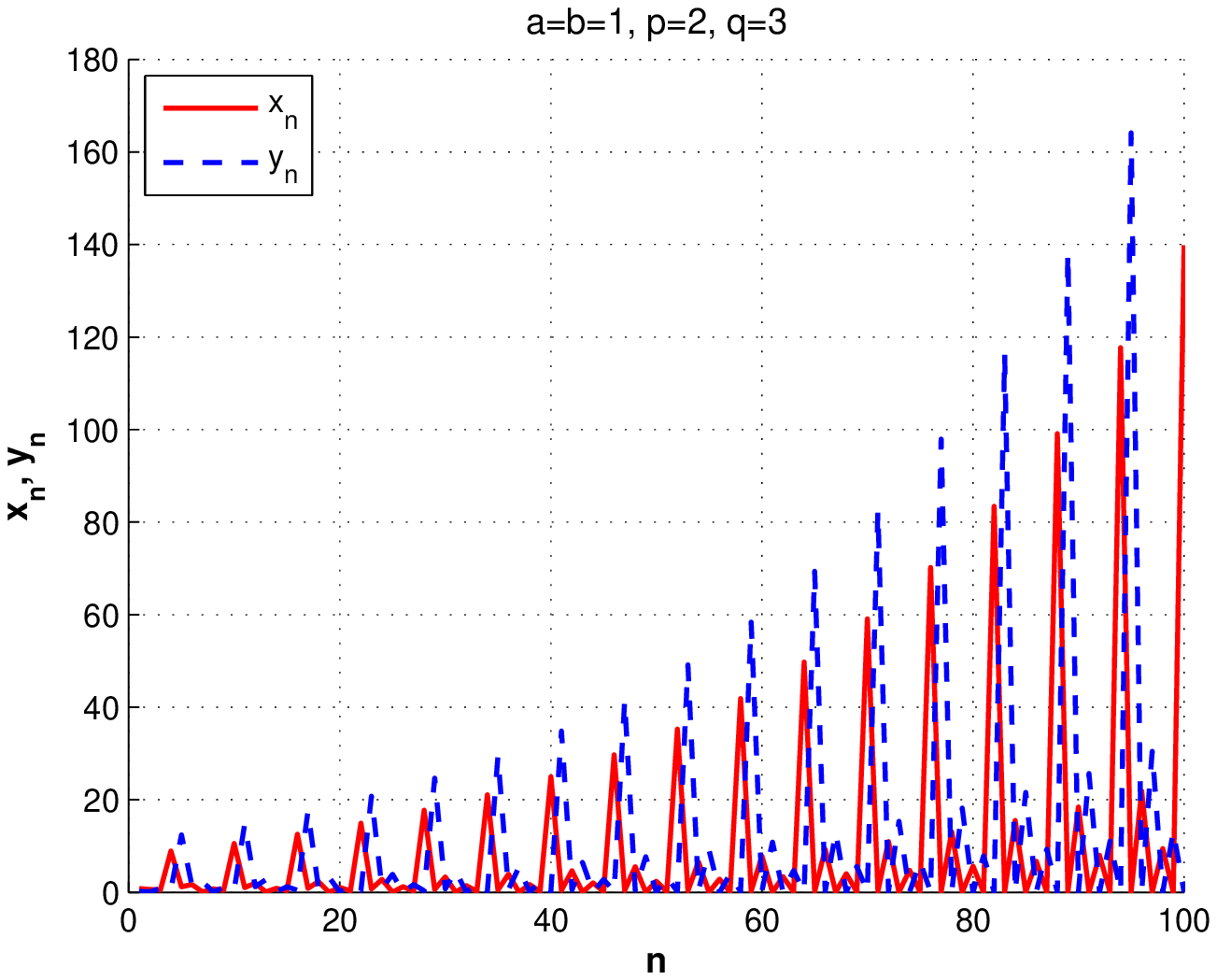}}\label{fig1}
     \caption{It is noticeable from the above illustration that the system of difference equation
    \[
    x_n=\frac{1}{y_{n-2}}, \quad y_n=\frac{y_{n-2}}{x_{n-3}y_{n-3}},\qquad n=0,1,\ldots
    \]
    has a solution $\{x_n\}$ (resp. $\{y_n\}$) with a subsequence $\{x_{6n+s}\}$ (resp. $\{y_{6n+t}\}$) ($0\leq s,t<6$) that tends to infinity.
    This result, moreover, agrees with Theorem \ref{main}--(i)}
\end{figure}

\begin{figure}
    \scalebox{.6}{\includegraphics{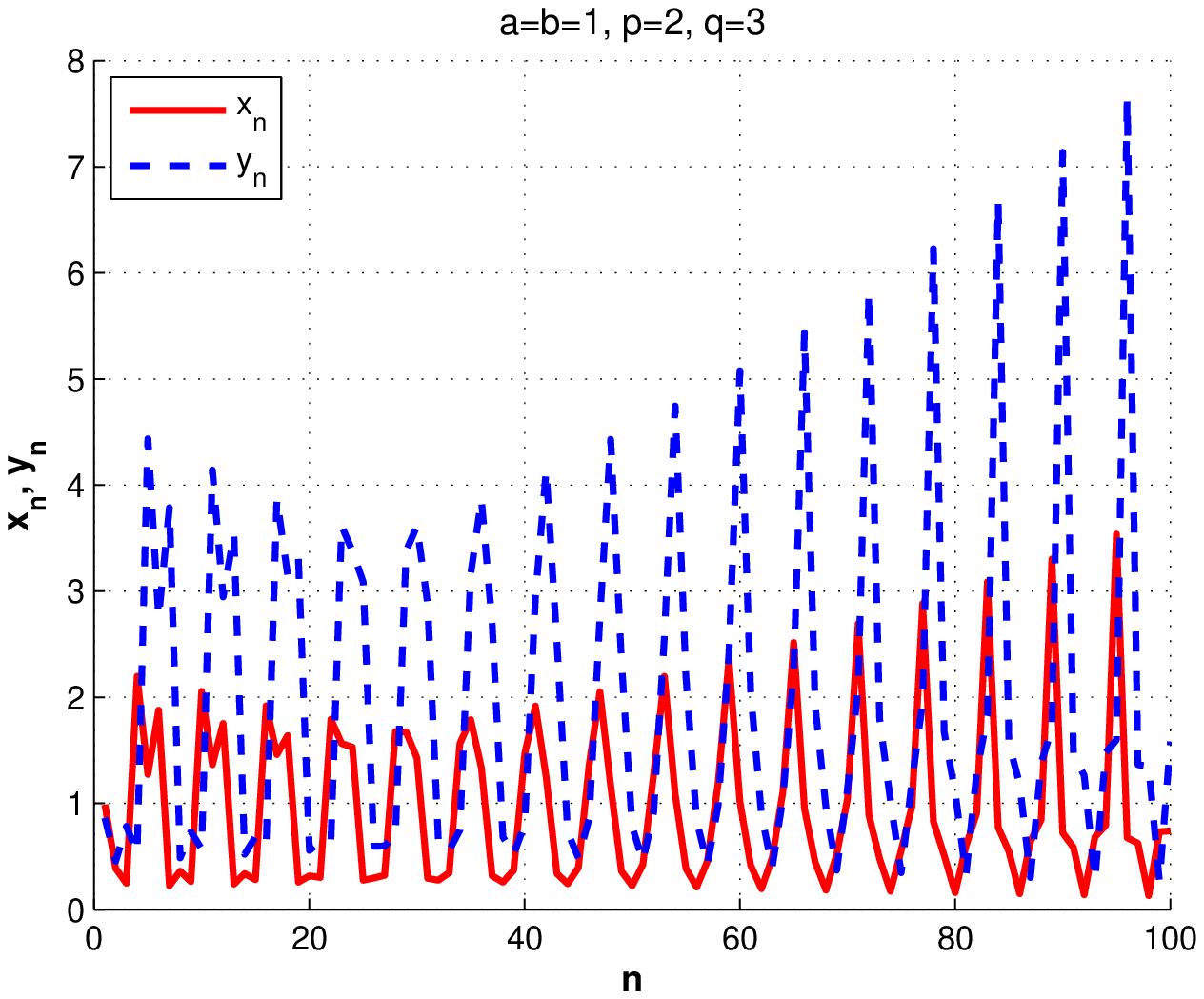}}\label{fig2}
    \caption{The above figure illustrates the behavior of another solution of the system
    \[
    x_n=\frac{1}{y_{n-2}}, \quad y_n=\frac{y_{n-2}}{x_{n-3}y_{n-3}},\qquad n=0,1,\ldots.
    \]
    Observe that after $N=25$, there is some subsequences $\{x_{6n+s}\}$ and $\{y_{6n+t}\}$ which both tends to infinity for $n \geq N$}
\end{figure}

\begin{figure}
    \scalebox{.6}{\includegraphics{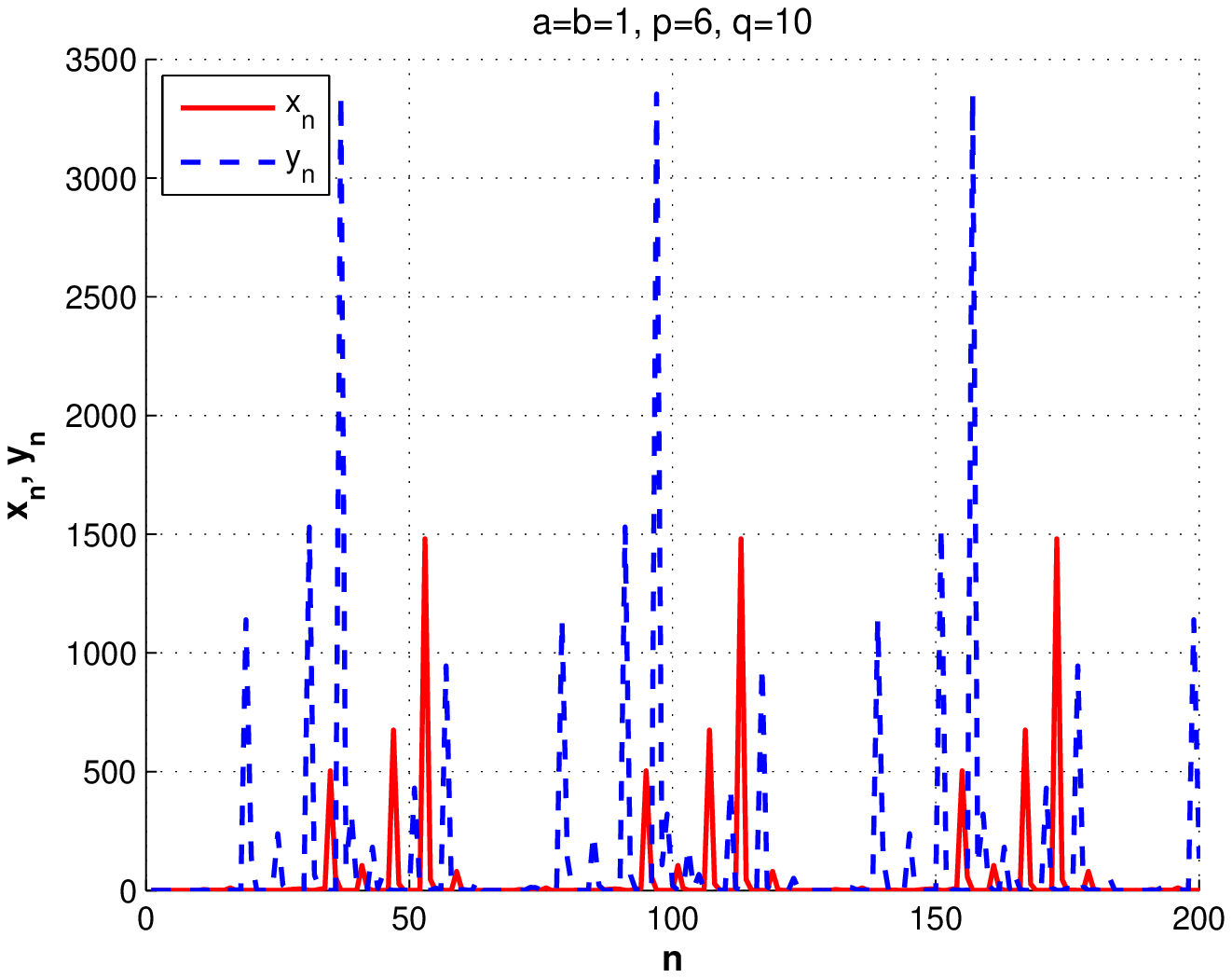}}\label{fig3}
    \caption{The above figure shows an interesting behavior of solutions of the system
    \[
    x_n=\frac{1}{y_{n-6}}, \quad y_n=\frac{y_{n-6}}{x_{n-10}y_{n-10}},\qquad n=0,1,\ldots.
    \]
    Clearly, as the above system satisfies the conditions in Theorem \ref{main}--(ii), we then have a periodic solution of period ${\rm lcm}(6,2\times 10)=60$}
\end{figure}

\begin{figure}
    \scalebox{.6}{\includegraphics{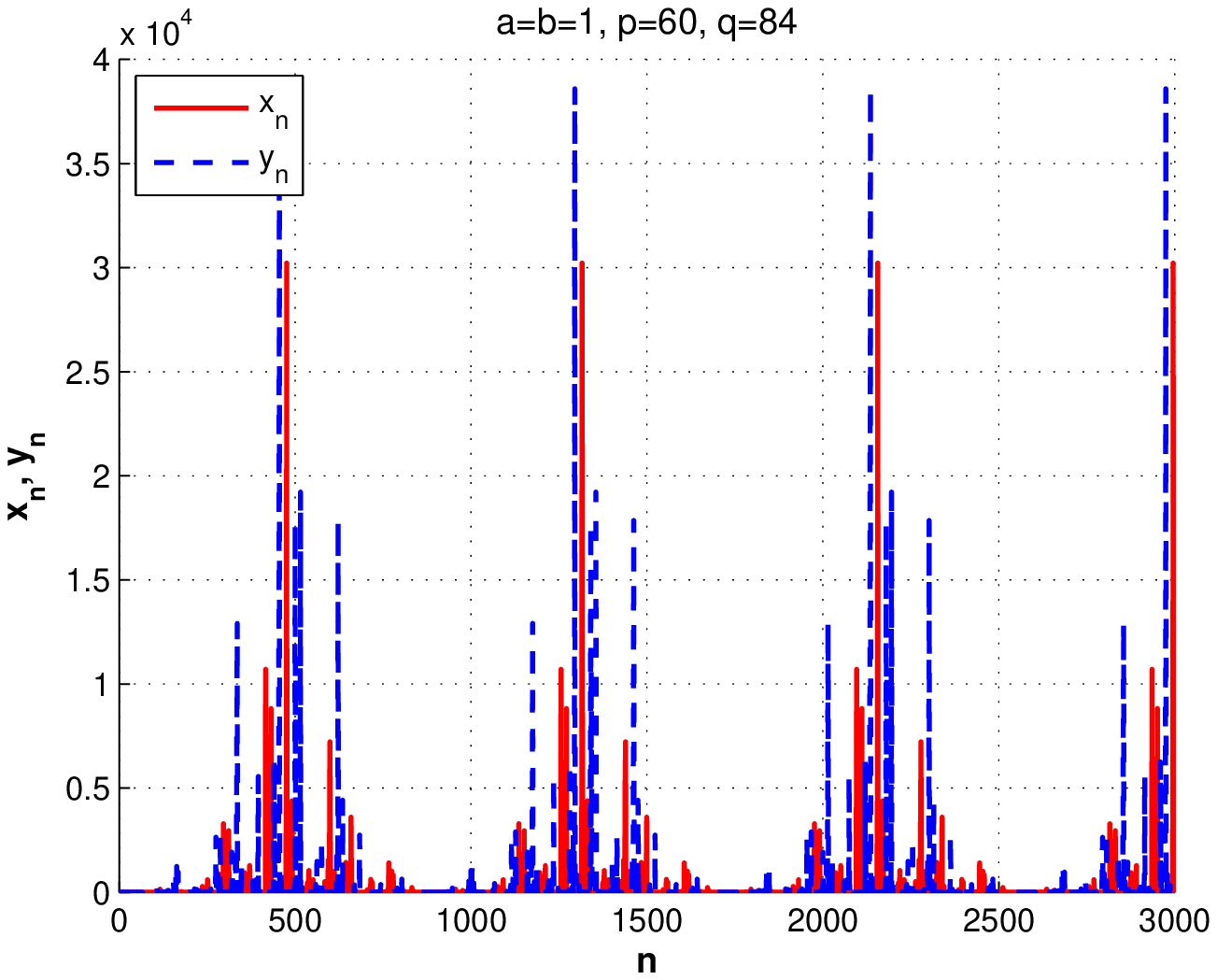}}\label{fig4}
    \caption{Another interesting illustration for Theorem \ref{main}--(ii) is shown above.
    In this example, the system
    \[
    x_n=\frac{1}{y_{n-60}}, \quad y_n=\frac{y_{n-60}}{x_{n-84}y_{n-84}},\qquad n=0,1,\ldots,
    \]
    has been considered. A simple computation for the period $m$ of the solution gives us $m={\rm lcm}(60,2\times 84)=840$}
\end{figure}

\begin{figure}
    \scalebox{.6}{\includegraphics{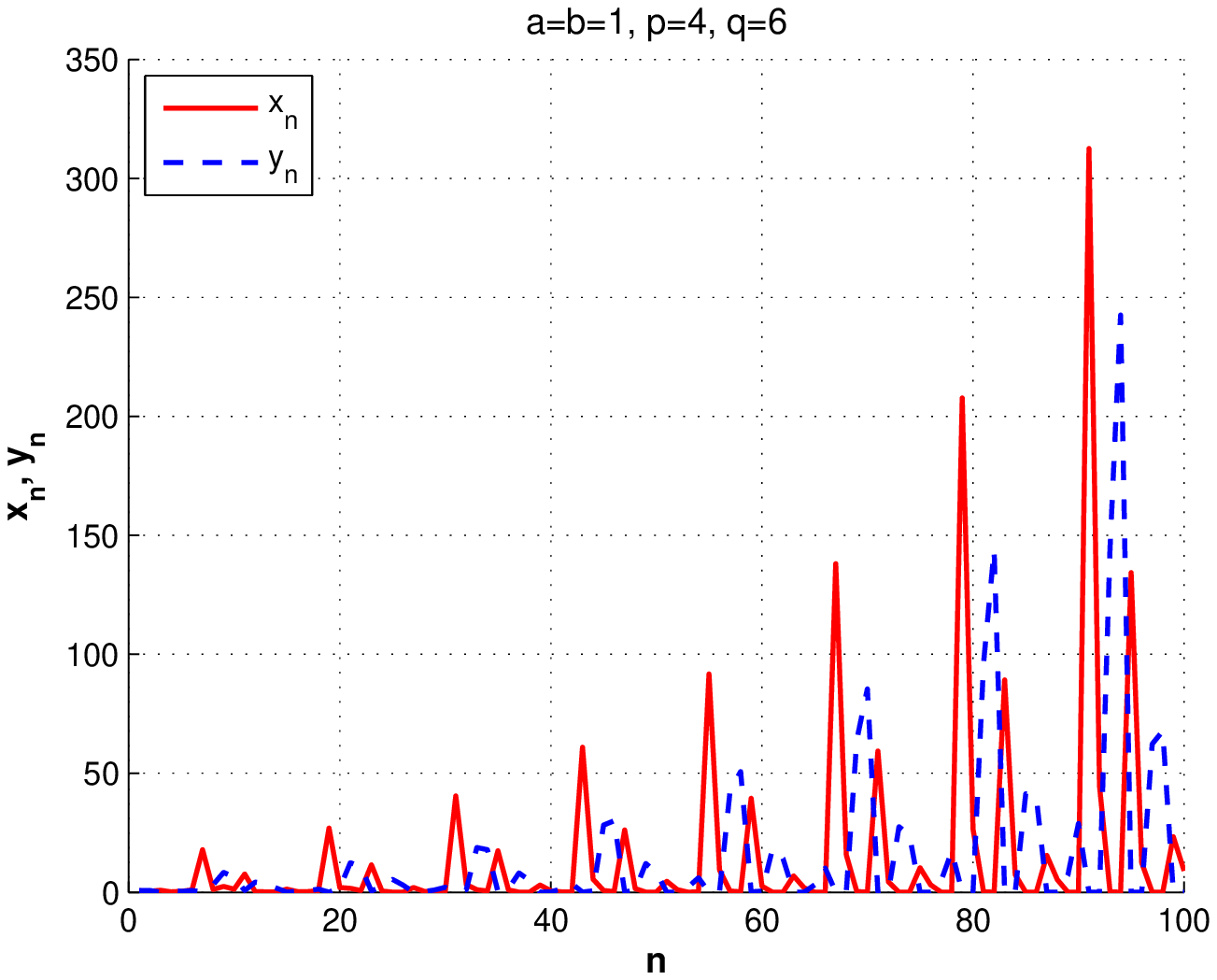}}\label{fig5}
    \caption{The above figure illustrates the behavior of a particular solution of the system
    \[
    x_n=\frac{1}{y_{n-4}}, \quad y_n=\frac{y_{n-4}}{x_{n-6}y_{n-6}},\qquad n=0,1,\ldots.
    \]
    Evidently, that particular solution of the given system behaves accordingly to Theorem \ref{main}--(iii)}
\end{figure}

\begin{remark}
We remark that a similar system has been studied by \"{O}zban in \cite{oz2} wherein he investigated the behavior of positive solutions of the system
	\[
		x_{n+1}=\frac{1}{y_{n-k}},\qquad y_{n+1}=\frac{y_{n-k}}{x_{n-m}y_{n-m-k}},\qquad n=0,1,\ldots,
	  \]
where $k$ is a nonnegative integer and $m$ is a positive integer. 
His main result states that all solutions of the above system of difference equations are periodic with period $2m+ 2k +2$. 
In particular, when $b=a$ and $k=0$, all solutions of above equation is periodic with period $2q+2$.
\end{remark}

\section{Conclusion}
In this short note, we have investigated the behavior of positive solutions of the system 
\[
	x_n=\frac{a}{y_{n-p}},\qquad y_n=\frac{by_{n-p}}{x_{n-q}y_{n-q}},\qquad n=0,1,\ldots,
\]    
for $b=a$, where $q$ is a positive integer with $p < q$, $p \nmid q$, and $p$ is an even number. 
We have found that every solution of the above system when $b=a$, with $p>0$, is non-periodic and has a subsequence that grows/decays exponentially whenever $q$ is odd.
However, a periodic solution of the given system occurs when $\gcd(p,q)>1$ and $p/\gcd(p,q)$ is odd. 
In this case, the period of the solution appears to be equal to the least common multiple of $p$ and $2q$.
On the other hand, a similar behavior as for the case when $q$ is odd was observed when $\gcd(p,q)>1$ and $p/\gcd(p,q)$ is even. 
Consequently, our result settled the question raised by Yang et al. in \cite{ycs} about the behavior of solution of the given equation on the case when $p$ is even and $q>p$ in the given system. 

\section*{Appendix}
From the polynomial equation $P(\lambda)=(\lambda^p-1)(\lambda^q+1)=0$, it is clear that $\lambda=1$ is a simple root.
Hence, a particular solution to the non-homogeneous equation
\begin{equation}
\label{receq}
a_n+a_{n-q}-a_{n-p}-a_{n-p-q}=\ln c, \qquad c:=a/b, 
\end{equation}
has the form 
\[
a_n^p=An
\]
from which, by simple computation, leads to $A=\ln (c/2p)$.
Thus, if $\gcd(p,q)>1$ and $p/\gcd(p,q)$ is odd, then the general solution of equation \eqref{receq} takes the form
\begin{align*}
x_n = e^{a_n}= c^{n/2p} \exp &\left\{ \sum_{l=0}^{p-1} \left( \alpha_{l,1} \cos \frac{2l\pi n}{p} + \alpha_{l,2} \sin \frac{2l\pi n}{p}\right)\right.\nonumber\\
			&\qquad+ \left. \sum_{k=0}^{q-1} \left( \beta_{k,1} \cos \frac{(2k+1)\pi n}{q} + \beta_{k,2} \sin \frac{(2k+1)\pi n}{q}\right)\right\}.\label{sol}
\end{align*}
We can write $x_n$ as $x_n=c^{n/2p} \hat{x}_n$, where $\hat{x}_n$ denotes the positive solution of equation \eqref{receq} with $c=1$. 

From above discussion, together with Theorem \ref{main}--(ii), we get the following results.
\begin{theorem}
Assume that $c \in(0,1)$, $\gcd(p,q)>1$ and $p/\gcd(p,q)$ is odd and let $m=\gcd(p,2q)$, then every positive solution of \eqref{receq} converges geometrically to zero. 
Moreover, for each $ t \in \{0, 1, . . . , m - 1\}$, the subsequence $\{x_{mn+t}\}_{n\in\mathbb{N}_0}$ converges monotonically to zero as $n$ tends to infinity.
\end{theorem}

\begin{theorem}
Assume that $c >1$, $\gcd(p,q)>1$ and $p/\gcd(p,q)$ is odd and let $m=\gcd(p,2q)$, then every positive solution of \eqref{receq} tends to infinity. 
Moreover, for each $ t \in \{0, 1, . . . , m - 1\}$, the subsequence $\{x_{mn+t}\}_{n\in\mathbb{N}_0}$ grows exponentially as $n$ tends to infinity.
\end{theorem}


\end{document}